%% file: draft3.tex
\documentclass[final]{siamart1116}

\input{draft3-shared}

\ifpdf
\hypersetup{
  pdftitle={\TheTitle},
  pdfauthor={\TheAuthors}
}
\fi

\begin{document}

\maketitle

\begin{abstract}
	\input{abstract}
\end{abstract}



\input{introduction}
\input{related_work}

\input{wspai}
\input{mmf2}

\input{experiments}
\input{conclusion}

\bibliographystyle{siamplain}
\bibliography{references}
\end{document}

%% file: draft3-shared.tex
\usepackage{lipsum}
\usepackage{amsfonts}
\usepackage{amsmath}
\usepackage{graphicx}
\usepackage{epstopdf}
\usepackage{algorithmic}
\usepackage{color}
\usepackage{rbase}
\usepackage{paralist}
\usepackage{tikz}
\usepackage{cleveref}
\usepackage{ulem}
\usepackage{stmaryrd}
\ifpdf
  \DeclareGraphicsExtensions{.eps,.pdf,.png,.jpg}
\else
  \DeclareGraphicsExtensions{.eps}
\fi

\numberwithin{theorem}{section}

\newcommand{\TheTitle}{A generic multiresolution preconditioner for sparse symmetric systems}
\newcommand{\TheAuthors}{Pramod Kaushik Mudrakarta, Risi Kondor}

\headers{Generic multiresolution preconditioner for sparse symm.\;systems}{\TheAuthors}

\title{{\TheTitle}}

\author{
  Pramod Kaushik Mudrakarta\footnotemark[2] \and Risi Kondor\thanks{Department of Computer Science, The University of Chicago, Chicago, IL}
}

\usepackage{amsopn}

\newcommand{\inv}{^{-1}}

\newcommand{\R}{\mathbb{R}}

\newcommand{\norm}[1]{\left|\left|#1\right|\right|}

\newcommand{\bl}[1]{\llbracket #1 \rrbracket} 
\newcommand{\offdiag}{\text{off-diag}}

\newcommand{\cupdot}{\mathbin{\mathaccent\cdot\cup}}
\newsiamthm{notation}{Notation}

\newcommand{\tikmxA}[1]{
	\br{\,\begin{tikzpicture}[baseline=-15.5, scale=0.045]
		\filldraw[gray] (0,0) rectangle +(#1,-#1); \end{tikzpicture}\,}}
\newcommand{\tikmxB}[2]{
	\br{\,\begin{tikzpicture}[baseline=-15.5, scale=0.045] 
		\filldraw[gray] (0,0) rectangle +(#1,-#1);
		\foreach \i in {#1,...,#2}{
			\filldraw[gray] (\i,-\i) rectangle +(1,-1);}
		\end{tikzpicture}\,}}


%% file: abstract.tex
We introduce a new general purpose multiresolution preconditioner for symmetric linear systems. 
Most existing multiresolution preconditioners use some standard wavelet basis that 
relies on knowledge of the geometry of the underlying domain. In constrast, based on the recently proposed Multiresolution Matrix Factorization (MMF) algorithm \cite{kondor2014multiresolution}, we construct a preconditioner that discovers a custom wavelet basis adapted to the given linear system 
without making any geometric assumptions.  Some advantages of the new approach are fast preconditioner-vector products, invariance to the ordering of the rows/columns, and the ability to handle systems of any size. Numerical experiments on finite difference discretizations of model PDEs and 
off-the-shelf matrices illustrate the effectiveness of the MMF preconditioner.

%% file: introduction.tex
\section{Introduction}
\label{sec:introduction}
Symmetric linear systems of the form 
\begin{align}
\label{eq:Ax=b} 
Ax = b, 
\end{align}
where $A\in \mathbb{R}^{n\times n}$ and $b\in \mathbb{R}^n$ are central to many numerical computations in 
science and engineering. Examples include finite difference discretizations of partial differential equations \cite{quarteroni2008numerical} and optimization algorithms where a linear system is solved in each 
iteration \cite{ye1991interior,ye1992finite,todd1990centered}.
Often, solving the linear system is the most time consuming part of large scale computations.

When $A$, the coefficient matrix, is large and sparse, usually iterative algorithms such as the minimum residual method (MINRES) \cite{paige1975solution} 
or the stabilized bi-conjugate gradient method (BiCGStab) \cite{van1992bi} are used to solve \cref{eq:Ax=b}.  
However, if the condition number $\kappa_2(A)$ is high (i.e., $A$ is ill-conditioned), 
these methods tend to converge slowly. For example, in the case of MINRES, for 
positive definite $A$, 
\begin{align}
\norm{Ax_n-b}_2 \le \left(1 - \frac{1}{\kappa_2(A)^2}\right)^\frac{n}{2}\norm{Ax_0-b}_2, 
\end{align}
where $x_n$ is the $n$-th iterate and $x_0$ is the initial guess \cite{trefethen1997numerical}.
Many matrices arising from problems of interest are ill-conditioned.  

Preconditioning is a technique to improve convergence, where, instead of \cref{eq:Ax=b}, we solve 
\begin{align}\label{eq:MAx=Mb}
 M\nts Ax = Mb,
\end{align}
where $M\in\mathbb{R}^{n\times n}$ is a rough approximation to $A\inv$ \footnote{An alternate way to precondition is from the right, i.e., solve $AMx = b$, but, for simplicity, in this paper we constrain ourselves to discussing left preconditioning.}. 
While \cref{eq:MAx=Mb} is still a large linear system, it is generally easier to solve than \cref{eq:Ax=b}, 
because $M\nts A$ is more favorably conditioned than $A$. 
Note that solving \cref{eq:MAx=Mb} with an iterative method involves computing 
many matrix-vector products with $M\nts A$, but that does not necessarily mean that $M\nts A$ needs to be computed explicitly.  
This is an important point, because even if $A$ is sparse, $M\nts A$ can be dense, and therefore expensive to compute. 

There is no such thing as a ``universal'' preconditioner. 
Preconditioners are usually custom-made for different kinds of coefficient matrices and are evaluated 
differently based on what kind of problem they are used to solve  
(how accurate $x$ needs to be, how easy the solver is to implement  
on parallel computers, storage requirements, etc.). 
Some of the most effective preconditioners exploit sparsity. 
The best case scenario is when both $A$ and $M$ are sparse, since in that case all matrix-vector 
products involved in solving \cref{eq:MAx=Mb} can be evaluated very fast. 
Starting in the 1970s, this lead to the devevelopment of so-called Sparse Approximate Inverse (SPAI) 
preconditioners \cite{benson1973iterative,grote1997parallel,benzi1999comparative,hawkins2005implicit}, 
which formuate finding $M$ as a  least squares problem 
\begin{align}
\label{eq:AM_I}
\min_{M\in\mathcal{S}}||AM-I||_F,
\end{align}
where $\mathcal{S}$ is an appropriate class of sparse matrices. 
Note that  since $||AM-I||^2_F = \sum_{i=1}^{n}||Am_i -e_i||^2_2$,
where $m_i$ is the $i$-th column of $M$ and $e_i$ is the $i$-th standard basis vector, 
\cref{eq:AM_I} reduces to solving $n$ independent least square problems, which can be done in parallel. 

One step beyond generic SPAI preconditioners are methods that use prior knowledge about the system at hand 
to transform $A$ to a basis where its inverse can be approximated in sparse form. 
For many problems, orthogonal wavelet bases are a natural choice. 
Recall that wavelets are similar to Fourier basis functions, but have the advantage of 
being localized in space. 
Transforming \rf{eq:Ax=b} to a wavelet basis amounts to rewriting it as 
$\widetilde{A}\tts \widetilde{x} = \widetilde{b}$, 
where 
\[\widetilde{A} = W^T\nts\nts A\tts W, \qquad\qquad \widetilde{x} = W^Tx, \qquad  \textrm{and}\qquad \widetilde{b} = W^T\tts b.\] 
Here, the wavelets appear as the columns of the orthogonal matrix $W$. 
This approach was first proposed by Chan, Tang and Wan \cite{chan1997wavelet}. 

Importantly, many wavelets admit fast transforms, meaning that $W^T$ factors in the form 
\begin{align}\label{eq:FWT}
W^T=W_L^T\tts W_{L-1}^T\ldots W_1^T, 
\end{align}
where each of the $W_\ell^T$ factors are sparse. While the wavelet transform itself 
is a dense transformation, in this case, 
transforming to the wavelet basis inside an interative solver 
can be done by sparse matrix-vector arithmetic exclusively. 
Each $W_\ell$ matrix can be seen as being responsible for extracting information from 
$x$ at a given scale, hence wavelet transforms constitute a form of multiresolution analysis.  


Wavelet sparse preconditioners have proved to be effective primarily in the PDE domain, where the 
problem is low dimensional and the 
structure of the equations (together with the discretization) strongly suggest the form of the wavelet transform. 
However, multiscale data is much more broadly prevalent, e.g., in biological problems and 
social networks. For these kinds of data, the underlying generative process is unknown, 
rendering the classical wavelet-based preconditioners ineffective. 

In this paper, we propose a preconditioner based on a form of multiresolution analysis for matrices 
called Multiresolution Matrix Factorization (MMF), that was first introduced in 
\cite{kondor2014multiresolution}. 
Similar to \cref{eq:FWT}, MMF has a corresponding fast wavelet transform, in particular, it is based on an 
approximate factorization of \m{A} of the form 
\begin{align}
A \approx Q_1^T\tts Q_2^T\ldots Q_L^T\tts H\tts Q_L\tts Q_{L-1}\ldots Q_1,
\end{align}
where each of the $Q_\ell$ matrices are sparse and orthogonal, and $H$ is close to diagonal. 
However, in contrast to classical wavelet transforms, here the $Q_\ell$ matrices are not induced 
from any specific analytical form of wavelets, but rather ``discovered'' by the algorithm itself 
from the structure of $A$, somewhat similarly to algebraic multigrid methods \cite{ruge1987algebraic}. 
This feature gives our preconditioner considerably more flexibility than existing 
wavelet sparse preconditioners, and allows it to exploit latent multiresolution structure in a wide 
range of problem domains. 

\subsection*{Notations}

In the following, we use \m{\sqb{n}} to denote the set \m{\cbrN{\oneton{n}}}. 
Given a matrix \m{A\tin\RR^{n\times n}} and two (ordered) sets \m{S_1,S_2\subseteq [n]},~\;  
\m{A_{S_1,S_2}} will denote the f\m{\absN{S_1}\<\times \absN{S_2}} dimensional 
submatrix of \m{A} cut out by the rows indexed by \m{S_1} and the columns indexed by \m{S_2}. 
\m{\wbar{S}} will denote the complement of \m{S}, in \m{[n]}, i.e., \m{[n]\!\setminus\! S}. 

%% file: related_work.tex
\section{Related work}
Constructing a good preconditioner hinges on two things: 
1. being able to design an efficient algorithm to compute an approximate inverse to $A$, and 
2. making the preconditioner as close to $A^{-1}$ as possible. 
It is rare for both a matrix and its inverse to be sparse. For example, Duff et al.\;\cite{duff1988sparsity} show that the inverses of irreducible, structurally sparse matrices are generally structurally dense. However, it is often the case that many entries of the inverse are small, making it possible to construct a good sparse approximate inverse. For example, \cite{demko1984decay} shows that when $A$ is banded and symmetric positive definite, the distribution of the magnitudes of the matrix entries in $A\inv$ decays exponentially. Benzi and Tuma \cite{benzi1999comparative} note that sparse approximate inverses have limited success because of the requirement that the actual inverse of the matrix has small entries.

A better way of computing approximate inverses is in factorized form using sparse factors. The dense nature of the inverse is still preserved in the approximation as the product of the factors (which is never explicitly computed) can be dense. Factorized approximate inverses have been proposed based on LU factorization. However, they are not easily parallelizable and are sensitive to reordering \cite{benzi1999comparative}. 

Multiscale variants of classic preconditioners have already been proposed and have often been found to be superior \cite{benzi2002preconditioning} to their one-level counterparts. 
The current frontiers of research on preconditioning also focus on designing algorithms for multi-core machines. Multilevel preconditioners assume that the coefficient matrix has a hierarchy in structure. These include the preconditioners that are based on rank structures, such as $\mathcal{H}$-matrices \cite{hackbusch1999sparse}, which represent a matrix in terms of a hierarchy of blocked submatrices where the off-diagonal blocks are low rank. This allows for fast inversion and LU factorization routines. Preconditioners based on $\mathcal{H}$-matrix approximations have been explored in \cite{faustmann2015mathcal,kriemann2015mathcal,grasedyck2009domain}. 
Other multilevel preconditioners based on low rank have been proposed in \cite{xi2016algebraic}.

Multigrid preconditioners \cite{briggs2000multigrid, pereira2006fast} are reduced tolerance multigrid solvers, which alternate between fine- and coarse-level representations to reduce the low and high frequency components of the error respectively. In contrast, hierarchical basis methods \cite{yserentant1986multi,yserentant1986hierarchical} precondition the original linear system as in \cref{eq:MAx=Mb} by expressing $A$ in a hierarchical representation. A hierarchical basis-multigrid preconditioner has been proposed in \cite{bank1988hierarchical}. 

Hierarchical basis preconditioners can be thought of as a special kind of wavelet preconditioners as it is possible to interpret the piecewise linear functions of the hierarchical basis as wavelets. Connections between wavelets and hierarchical basis methods have also been explored in \cite{vassilevski1997stabilizing,vassilevski1998stabilizing} to improve the performance of hierarchical basis methods.

%% file: wspai.tex
\section{Wavelet based sparse approximate inverse preconditioners}
\label{sec:wspai}

We begin with a brief introduction to classical orthogonal wavelet transforms. 
For a detailed introduction, see \cite{daubechies1992ten}. 
Assuming $n=2^N$ for simplicity,
the $L$-level wavelet transform of a signal $x\in\R^n$ \ignore{where $n=2^N$} can be written as a matrix 
vector product $W^Tx$, where
\begin{align}
W = W_1 W_2\ldots W_L
\end{align}
with $L\leq n$ and 
\begin{align}
W_k^T = \begin{pmatrix}
U_k & 0 \\
V_k & 0 \\
0 & I_{n - \frac{n}{2^k-1}}
\end{pmatrix}
\hspace{70pt} k=1,\ldots,L,  
\end{align}
where $U_k, V_k\in\R^{(n/2^k)\times (n/2^{k-1}) }$ are of the form
\begin{align*}
U_k &= \begin{pmatrix}
h_0 & h_1 & h_2    & \cdots & h_{m-1} &&& \\
&     & h_0    & h_1    & h_2     & \cdots & h_{m-1} & \\
&	  &	\ddots & \ddots & \ddots  & \ddots &         & \\
h_2 & \cdots & h_{m-1} &&&& h_0 & h_1       
\end{pmatrix},\\ 
&\\
V_k &= \begin{pmatrix}
g_0 & g_1 & g_2    & \cdots & g_{m-1} &&& \\
&     & g_0    & g_1    & g_2     & \cdots & g_{m-1} & \\
&	  &	\ddots & \ddots & \ddots  & \ddots &         & \\
g_2 & \cdots & g_{m-1} &&&& g_0 & g_1       
\end{pmatrix}.
\end{align*} 
The scalars $h_i, g_i$ for $i=1,\ldots,m-1$ are the high-pass and low-pass filter coefficients of the wavelet 
transform, respectively. The above holds true even when $n = p\ts 2^s$ for some $s$ and $p$. 
In that case, the maximum level of the wavelet transform applied is upper bounded by $s$. 

On higher dimensional signals, wavelet transforms are applied dimension-wise. 
For example, let $x\in\R^{n^2}$ be a 2D signal (matrix) which has been vectorized by stacking the columns. 
The wavelet transform $\widetilde{x}$ is computed by first applying a 1D transform on the columns and then on the rows. 
If $W\in\R^{n\times n}$ is the 1D orthogonal wavelet transform matrix, then 
\begin{align*}
\widetilde{x} = (I_n\otimes W^T)(W^T\otimes I_n)\, x = (W \otimes W)^T x, 
\end{align*} 
where $\otimes$ is the Kronecker product \cite{van2000ubiquitous} and $I_n$, the $n\<\times n$ identity matrix. 
Thus, $W\otimes W$ can be called the two dimensional wavelet transform matrix. 
For vectorized 3D signals (tensors), the wavelet transform matrix is $W\otimes W\otimes W$.

Chan, Tang and Wan \cite{chan1997wavelet} were the first to propose a wavelet sparse approximate inverse preconditioner. 
In their approach, the linear system \cref{eq:Ax=b} is first transformed into a standard wavelet basis such 
as Daubechies the \cite{daubechies1992ten} basis, 
and a sparse approximate inverse preconditioner is computed for the transformed coefficient matrix by solving 
\begin{align}
\label{eq:chan_tang_wan_minimization_problem}
\min_{M\in\mathcal{S}_\text{blockdiag}} \norm{\,WAW^TM - I\,}_F.
\end{align}
The preconditioner is constrained to be block diagonal in order to maintain its sparsity and simplify computation. They show the superiority of the wavelet preconditioner over an adaptive sparse approximate inverse preconditioner for elliptic PDEs with smooth coefficients over regular domains. 
However, their method performs poorly for elliptic PDEs with discontinuous coeffcients. The block diagonal constraint does not fully capture the structure of the inverse in the wavelet basis.

Bridson and Tang \cite{bridson2001multiresolution} construct a multiresolution preconditioner similar to Chan, Tang and Wan \cite{chan1997wavelet}, but determine the sparsity structure adaptively. Instead of using Daubechies wavelets, they use second generation wavelets \cite{sweldens1998lifting}, which allows the preconditioner to be effective for PDEs over irregular domains. However, their algorithm requires the additional difficult step of finding a suitable ordering of the rows/columns of the coefficient matrix which limits the number of levels to which multiresolution structure can be exploited. 

Hawkins and Chen \cite{hawkins2005implicit} compute an implicit wavelet sparse approximate inverse preconditioner, which removes the computational overhead of transforming the coefficient matrix to a wavelet basis. Instead of \cref{eq:chan_tang_wan_minimization_problem}, they solve 
\begin{align}
\label{eq:hawkins_chen_minimization_problem}
\min_{M\in\mathcal{S}_W} \norm{WAM - I}_F,
\end{align}
where $\mathcal{S}_W$ is the class of matrices which have the same sparsity structure as $W$. They empirically show that this sparsity constraint is enough to construct a preconditioner superior to that of Chan, Tang and Wan \cite{chan1997wavelet}. The complete algorithm is described in \cref{alg:hawkins1,alg:hawkins2}.

\input{alg-spai1.tex}

\input{alg-spai2.tex}

Hawkins and Chen \cite{hawkins2005implicit} apply their preconditioner on Poisson and elliptic PDEs in 1D, 2D and 3D. 
We found, by experiment, that it is critical to use a wavelet transform of the same dimension as the underlying 
PDE of the linear system for success of their preconditioner. 
On linear systems where the underlying data generator is unknown --- this happens, for example, when we are 
dealing with Laplacians of graphs --- their preconditioner is ineffective. 
Thus, there is a need for a wavelet sparse approximate inverse preconditioner which can mould itself to any 
kind of data, provided that it is reasonable to assume a multiresolution structure.

%% file: alg-spai1.tex
\begin{algorithm}[t]
	\caption{Solve $Ax=b$ using the implicit wavelet SPAI preconditioner \cite{hawkins2005implicit}}
	\begin{algorithmic}[1]
		\STATE Compute preconditioner $\widehat{M} = \arg\min_{M\in\mathcal{S}_W} ||AM - W ||_F$
		\STATE Solve $W^T A \widehat{M} y = W^T b$
		\RETURN $x = \widehat{M}y$
	\end{algorithmic}
	\label{alg:hawkins1}
\end{algorithm}

%% file: alg-spai2.tex
\begin{algorithm}[t]
	\caption{Compute preconditioner $\widehat{M} = \arg\min_{M\in\mathcal{S}_W} ||AM - W ||_F$ }
	\begin{algorithmic}[1]
		\FOR{$j=1,\ldots,n$}
		\STATE $S_j = $ indices of nonzero entries of $w_j$
		\STATE $T_j = $ indices of nonzero entries of $A(:,S_j)$
		\STATE Solve $z^* = \arg\min ||A(T_j,S_j)z - w_j(T_j)||_2$ by reduced QR-factorization 
		\STATE Set $\widehat{m}_j(T_j) = z^*$
		\ENDFOR
		
		\RETURN $\widehat{M}$
	\end{algorithmic}
	\label{alg:hawkins2}
\end{algorithm}

%% file: mmf2.tex
\section{Multiresolution matrix factorization}
\label{sec:mmf}
The Multiresolution Matrix Factorization (MMF) of a symmetric matrix $A\in\R^{n\times n}$, 
as defined in \cite{kondor2014multiresolution}, is a multilevel sparse factorization of the form
\begin{align}
\label{eq:mmf}
A \approx Q_1^T Q_2^T \ldots Q_L^T H\tts Q_L \ldots Q_2 Q_1, 
\end{align}
where the matrices $Q_1,\ldots,Q_L$ and $H$ obey the following conditions:
\begin{enumerate}[~1.]
\item Each $Q_\ell$ is orthogonal and highly sparse. 
In the simplest case, each $Q_\ell$ is a Givens rotation, i.e., a matrix which differs from the identity in 
just the four matrix elements
\begin{align*}
[Q_\ell]_{i,i} = \cos \theta,&\qqquad [Q_\ell]_{i,j} = - \sin \theta,\\
[Q_\ell]_{j,i} = \sin \theta,&\qqquad [Q_\ell]_{j,j} = \cos \theta,
\end{align*}
for some pair of indices $(i,j)$ and rotation angle $\theta$. 
Multiplying a vector with such a matrix rotates it counter-clockwise by $\theta$ in the $(i,j)$ plane. 
More generally, $Q_\ell$ is a so-called 
$k$-point rotation, 
which rotates not just two, but $k$ coordinates. 
\item Typically, in MMF factorizations $L\<=O(n)$, and the size of the active part of the \m{Q_\ell} matrices decreases according to a set schedule $n = \delta_0 \ge \delta_1 \ge \ldots \ge \delta_L$. 
More precisely, there is a nested sequence of sets $[n] = S_0 \supseteq S_1 \supseteq \ldots \supseteq S_L$ 
such that the $[Q_\ell]_{\overline{S_{\ell-1}},\overline{S_{\ell-1}}}$ 
part of each rotation is the $n-\delta_{\ell-1}$ 
dimensional identity. $S_\ell$ is called the active set at level $\ell$. In the simplest case,  $\delta_\ell=n-\ell$.
\item $H$ is an $S_L$-core-diagonal matrix, which means that it   
is block diagonal with two blocks: $H_{S_L,S_L}$, called the core, which is dense, and 
$H_{\wbar{S_L},\wbar{S_L}}$ which is diagonal. In other words, $H_{i,j}\!=\!0$ unless $i,j\tin S_L$ or $i\!=\!j$.  
\end{enumerate}
The structure implied by the above conditions is illustrated in \cref{fig:mmf}. MMF factorizations are, 
in general, only approximate, as there is no guarantee that 
$O(n)$ sparse orthogonal matrices can bring a symmetric matrix to core-diagonal form.
Rather, the goal of MMF algorithms is to minimize the approximation error, which, in the simplest case, 
is the Frobenius norm of the difference between the original matrix and its MMF factorized form. 
\input{fig-mmf}

MMF was originally introduced in the context of multiresolution analysis on discrete spaces, 
such as graphs. In particular, the columns of 
$Q^T=Q_1^T \ldots Q_{L-1}^T Q_L^T$ have a natural interpretation as wavelets, and the factorization 
itself is effectively a fast wavelet transform, mimicking the structure of classical orthogonal 
multiresolution analyses on the real line \cite{mallat1989theory}. 
MMF has also been successfully used for compressing large matrices \cite{teneva2016multiresolution}. 

In this paper we use MMF in a different way. The key property that we exploit is that 
\cref{eq:mmf} automatically gives rise to an approximation to $A\inv$,  
\begin{align}\label{eq:immf}
\widetilde{A^{-1}}=Q_1^T\ldots Q_{L-1}^T\tts Q_L^T\tts H\inv\tts Q_L\tts Q_{L-1}\ldots Q_1,  
\end{align}
which is very fast to compute, since inverting $H$ reduces to separately inverting its core 
(which is assumed to be small) and inverting its diagonal block (which is trivial).  
Assuming that the core is small enough, the overall cost of inversion becomes $O(n)$. 
When using \cref{eq:immf} as a preconditioner, of course we never compute \cref{eq:immf} 
explicitly, but rather (similarly to other wavelet sparse approximate inverse preconditioners)  
we apply it to vectors in factorized form as  
\begin{align}
\widetilde{A^{-1}}\ts v=Q_1^T(\ldots (Q_{L-1}^T (Q_L^T (H^{-1} (Q_L (Q_{L-1}\ldots (Q_1v))\ldots). 
\end{align}	
Since each of the factors here is sparse, the entire product can be computed in \m{O(n)} time.

\ignore{
MMF offers key computational advantages. While the product on the r.h.s of \cref{eq:mmf} may be dense, the rotation matrices $Q_1,\ldots,Q_L$ and the center matrix $H$ are highly sparse which allows for fast matrix-vector products. Product of the MMF with a vector $v$ is simply
\begin{align}
Q_1^T(\ldots (Q_{L-1}^T (Q_L^T (H (Q_L (Q_{L-1}\ldots (Q_1\tts v))\ldots),
\end{align}
The inverse of an MMF is 
\begin{align}
Q_1^T\ldots Q_{L-1}^T Q_L^T H\inv Q_L Q_{L-1}\ldots Q_1,
\end{align}
whose computational cost is essentially the cost of inverting of the dense part of the core-diagonal matrix $H$, which can be done via eigenvalue decomposition. As we restrict the size of the dense part of $H$ to be very small and inverting the diagonal part is simply computing the reciprocals of the diagonal elements, the inverse of an MMF can essentially be done in $O(n)$ time. Thus, the preservation of the factorization structure and low computational costs allow the use of the MMF inverse as a fast preconditioner for linear systems. 
}


\subsection*{Computation of the MMF}
\label{sec:computation_mmf}

The MMF of a symmetric matrix $A$ is usually computed by minimizing the Frobenius norm
factorization error  
\begin{align}
\label{eq:mmf_minimization_problem}
\nmN{\ts A - Q_1^T\ldots Q_L^T H\tts Q_L\ldots Q_1\ts }_{\text{Frob}}
\end{align}
over all admissible choices of active sets $S_1,\ldots,S_L$ and rotation matrices $Q_1,\ldots,Q_L$. 
The minimization is carried out in a greedy manner, 
where the rotation matrices $Q_1,\ldots,Q_L$ are determined sequentially, as $A$ is 
subjected to the sequence of transformations 
\[A\mapsto \underbrace{Q_1 A\ts Q_1^T}_{A_1} \mapsto \underbrace{Q_2 Q_1 A\ts Q_1^T Q_2^T}_{A_2} \mapsto \ldots 
\mapsto \underbrace{Q_L\ldots Q_2 Q_1 A\ts Q_1^T Q_2^T\ldots Q_L^T}_{H}.\]  
In this process, at each level $\ell$, the algorithm 
\begin{enumerate}[~1.]
\item Determines which subset of rows/columns $\{i_1,\ldots,i_k\}\subseteq S_{\ell-1}$ 
are to be involved in the next rotation, $Q_\ell$. 
\item Given $\{i_1,\ldots,i_k\}$, it optimizies the actual entries of $Q_\ell$. 
\item Selects a subset of the indices in $\{i_1,\ldots,i_k\}$ for removal from the active set 
(the corresponding rows/columns of the working matrix $A_\ell$ 
then become ``wavelets'').  
\item Sets the off-diagonal parts of the resulting wavelet rows/columns to zero in $H$. 
\end{enumerate}
The final error is the sum of the squares of the zeroed out off-diagonal elements (see Proposition 1 in \cite{kondor2014multiresolution}). 
The objective therefore is to craft $Q_\ell$ such that these off-diagonals are as small as possible.

For preconditioning it is critical to be able to compute the MMF approximation fast. 
To this end employ two further heuristics. 
First, the row/column selection process is accelerated by randomization: for each $\ell$, the first index 
$i_1$ is chosen uniformly at random from the current active set $S_{\ell-1}$, and then 
$i_2,\ldots,i_k$ are chosen so as to ensure that $Q_\ell$ can produce $\delta_\ell-\delta_{\ell-1}$  
rows/columns with suitably small off-diagonal norm. 
Second, exploiting the fundamentally local character of MMF pivoting, the entire algorithm is 
parallelized using a generalized blocking strategy first described in \cite{teneva2016multiresolution}. 
 	
\begin{notation} 
Let 
$B_1 \cupdot B_2 \cupdot \ldots \cupdot B_k = [n]$ be a partition of $[n]$ and $A\in\RR^{n\times n}$. 
We use $\bl{A}_{i,j}$ to denote the $[A]_{B_i,B_j}$ block of $A$ and say that $A$ 
is $(B_1,\ldots,B_k)$-block-diagonal if $\bl{A}_{i,j}\<=0$ if $i\<\neq j$. 
\end{notation}

The pMMF algorithm proposed in \cite{teneva2016multiresolution} 
uses a rough clustering algorithm to group the rows/columns of $A$ into a certain number of blocks,   
and factors each block independently and in parallel. However, to avoid overcommitting to a specific 
clustering, each of these factorizations is only partial (typically the core size is on the order of 
$1/2$ of the size of the block). The algorithm proceeeds in stages, where each stage consists 
of (re-)clustering the remaining active part of the matrix, performing partial MMF on each 
cluster in parallel, and then reassembling the active rows/columns from each cluster into a single matrix again 
(\cref{alg:pMMF_top_level}). 
\input{alg-pmmf.tex}

Assuming that there are $P$ stages in total, this process results in a two-level factorization. 
At the stage level, we have 
\begin{align}
\label{eq:pMMF_expansion}
A \approx \overline{Q}_1^T\ts \overline{Q}_2^T\ldots\overline{Q}_P^T\ts H\ts \overline{Q}_P\ldots \overline{Q}_2\ts \overline{Q}_1,
\end{align}
where, assuming that the clustering in stage $p$ is \m{B^p_1\cupdot B^p_2\cupdot \ldots \cupdot B^p_m}, 
each $\overline{Q}_p$ is a \m{(B^p_1,\ldots,B^p_m)} block diagonal orthogonal matrix, which, in turn, 
factors into a product of a large number of elementary $k$-point rotations 
\begin{align}
\label{eq:pMMF_stage_expansion}
\overline{Q}_p = Q_{l_p}\ldots Q_{l_{p-1}+2}Q_{l_{p-1}+1}. 
\end{align}
Thanks to the combination of these computational tricks, 
empirically, for sparse matrices, pMMF can achieve close to linear scaling behavior with $n$, both in 
memory and computation time \cite{teneva2016multiresolution}. 
For completeness, the subroutine used to compute the rotations in each cluster is presented in \cref{alg:pMMF_stage_comp}. 

\input{alg-findrot}

\ignore{
We consider the pMMF algorithm for minimizing \cref{eq:mmf_minimization_problem} 
which was first proposed by Teneva et al.\cite{teneva2016multiresolution}. It is a parallelized algorithm which computes an MMF in almost linear time with respect to the size of $A$. 
pMMF factors $A$ as
\begin{align}
\label{eq:pMMF_expansion}
A \approx \overline{Q}_1^T\overline{Q}_2^T\ldots\overline{Q}_P^T H \overline{Q}_P\ldots \overline{Q}_2 \overline{Q}_1
\end{align}
where each $\overline{Q}_p$ is a compound rotation which itself takes the form
\begin{align}
\label{eq:pMMF_stage_expansion}
\overline{Q}_p = Q_{l_p}\ldots Q_{l_{p-1}+2}Q_{l_{p-1}+1}
\end{align}
of a product of (a typically large number of) Givens rotations. 

To enable parallelism, all the Givens rotation matrices in \cref{eq:pMMF_stage_expansion} are forced to conform to the same block structure.

pMMF exploits parallelization and radically reduces the complexity of MMF computation. This is achieved by forcing $\overline{Q}_p$ to conform to a block-diagonal structure such that each $Q_i$ in \cref{eq:pMMF_stage_expansion} also has an identical block-diagonal structure. The computation of the blocks comprising $\overline{Q}_p$ are done independently and can be distributed over the available processors in the computer. Expanding each stage, \cref{eq:pMMF_expansion} is identical to \cref{eq:MMF_expansion} except that $Q_1,\ldots,Q_L$ are forced into contiguous runs of rotations having the same block structure.

At each stage, a blocked configuration of $A$ is available via clustering the columns based on normalized inner product. This ensures that columns of $A$ having high inner product with each other fall into the same cluster. The gram matrix $B^TB$ of each cluster $B$ is computed and columns having high inner product with each other are selected. A Givens rotation is constructed using these columns and applied to both the gram matrix and $A$. This process is repeated until a maximum number of rotations for that stage is reached. 

The pseudocode for the entire pMMF algorithm is described in \cref{alg:pMMF_top_level,alg:pMMF_stage_comp}
}

%% file: fig-mmf.tex
\begin{figure*}[t]
\label{fig:mmf}
\begin{equation*}
PAP^T\;\approx\;
\underset{\displaystyle Q_1^\top}{\tikmxA{17}}\nts 
\underset{\displaystyle Q_2^\top}{\tikmxB{14}{17}}\hdots 
\underset{\displaystyle Q_L^\top}{\tikmxB{5}{17}}
\underset{\displaystyle H^{\phantom{\top}}}{\tikmxB{3}{17}}\vspace{-15pt}
\underset{\displaystyle Q_L^{\phantom{\top}}}{\tikmxB{5}{17}}
\hdots \underset{\displaystyle Q_2^{\phantom{\top}}}{\tikmxB{14}{17}}\nts
\underset{\displaystyle Q_1^{\phantom{\top}}}{\tikmxA{17}}
\end{equation*}\vspace{-10pt}\mbox{}\\
\caption{\label{fig: mmf} 
A graphical representation of the structure of Multiresolution Matrix Factorization. 
Here, \m{P} is a permutation matrix which ensures that \m{S_\ell\<=\cbrN{1,\ldots,\delta_\ell}} for each \m{\ell}. 
Note that \m{P} is introduced only for the sake of visualization, 
an actual MMF would not contain such an explicit permutation. 
}\vspace{-12pt}
\end{figure*}

%% file: alg-pmmf.tex
\begin{algorithm}[t]
\begin{algorithmic}
\STATE \textbf{Input:}~ a symmetric matrix \m{A\tin\R^{n\times n}}
\STATE \m{A_0\leftarrow A}
\STATE \textbf{for} (\m{p\<=1} to \m{P})\,\m{\{}
\STATE ~~~~Cluster the active columns of \m{A_{p-1}} to \m{B^p_1\cupdot B^p_2 \cupdot \ldots \cupdot B^p_m}
\STATE ~~~~Reblock \m{\, A_{p-1}} according to \m{(B^p_1,\ldots,B^p_m)}
\STATE ~~~~\textbf{for} (\m{u\<=1}~\textbf{to}~\m{m})~~\m{\bl{\wbar{Q}_p}_{u,u}\leftarrow} 
\textsc{FindRotationsForCluster}(\m{[A_p]_{\,:\,,B_u}})
\STATE ~~~~\textbf{for}\,\m(\m{u\<=1}~\textbf{to}~\m{m})~\{
\STATE ~~~~~~~~\textbf{for}\,(\m{v\<=1}~\textbf{to}~\m{m})~\{
\STATE ~~~~~~~~~~~~\m{\bl{A_{p}}_{u,v}\nts\leftarrow {\bl{\wbar{Q}_p}_{u,u}}\bl{A_{p-1}}_{u,v} \bl{\wbar{Q}_p}_{v,v}{\!\!\!\!\!}^{\top}} 
\STATE ~~~~\}\}
\STATE \}
\STATE \m{H\leftarrow} the core of \m{A_L} plus its diagonal 
\STATE \textbf{Output:}~ \m{(H,\wbar{Q}_1,\ldots,\wbar{Q}_p)}
\end{algorithmic}
\caption{\label{alg:pMMF_top_level}~\textbf{pMMF}~\;(top level of the pMMF algorithm)}
\end{algorithm}

%% file: alg-findrot.tex
\begin{algorithm}[t]
\begin{algorithmic}
	\STATE \textbf{Input:}~ a matrix \m{\Ucal} made up of 
	the \m{c} columns of \m{A_{p-1}} forming cluster \m{u} in \m{A_p}\\ 
	\STATE Compute the Gram matrix \m{G\<=\Ucal^\top \Ucal}\\
	\STATE \m{S\leftarrow \cbrN{1,2,\ldots,c}}~~ (the active set) \\
	\STATE \textbf{for}~(\m{s\<=1}~\textbf{to}~\m{\lfloor \eta c\rfloor})\m{\{}\\
	\STATE ~~~~Select\; \m{i\tin S} uniformly at random\\ 
	\STATE ~~~~Find\; \m{j\<=\argmax_{S\setminus\cbr{i}} 
		\abs{\inpN{\Ucal_{:,i},\Ucal_{:,j}}}/\nmN{\ts \Ucal_{:,j}\ts}}\\ 
	\STATE ~~~~Find\; the optimal Givens rotation~\m{q_s}~of columns \m{(i,j)}\\ 
	\STATE ~~~~\m{\Ucal\leftarrow q_s\ts \Ucal\ts q_s^\top}\\
	\STATE ~~~~\m{G\leftarrow q_s\tts G\ts  q_s^\top}\\
	\STATE ~~~~\textbf{if}~\:\m{\nm{\ts \Ucal_{i,:}\ts }_{\offdiag}\! < \nm{\ts \Ucal_{j,:}\ts}_{\offdiag}}~~\textbf{then}~~\m{S\leftarrow S\setminus\cbr{i}}
	~~\textbf{else}~~\m{S\leftarrow S\setminus\cbr{j}}\\
	\STATE \m{\}}
	\STATE \textbf{Output:}~ \m{\bl{\wbar{Q}_p}_{u,u}=q_{\lfloor \eta c\rfloor}\ldots q_2\tts q_1}
\end{algorithmic}
\caption{~\;\textsc{FindRotationsForCluster}(\m{\Ucal)}~ (we assume \m{k\<=2} and 
	\m{\eta} is the compression ratio) 
\label{alg:pMMF_stage_comp}}
\end{algorithm}

%% file: experiments.tex
\section{Numerical results}
\label{sec:experiments}

We consider both model PDE problems and off-the-shelf datasets for comparing the preconditioners. The model PDE problems used are
\begin{itemize}
	\item \textit{1D Laplacian.} One dimensional Poisson's equation 
	\begin{align*}
	u_{xx} = (1+x^2)^{-1} e^x ,\quad x\in[0,1],
	\end{align*}
	with a Dirichlet boundary condition discretized with central differences.
	
	\item \textit{2D Laplacian.} Two dimensional Poisson's equation
	\begin{align*}
	u_{xx} + u_{yy} = -100x^2,\quad (x,y)\in [0,1]^2,
	\end{align*}
	with a Dirichlet boundary condition discretized with central differences.
	
	\item \textit{3D Laplacian.} Three dimensional Poisson's equation
	\begin{align*}
	u_{xx} + u_{yy} + u_{zz} = -100x^2,\quad (x,y,z)\in [0,1]^3.	
	\end{align*}
	
	\item \textit{2D Disc.} Two dimensional PDE with discontinuous coefficients
	\begin{align*}
	(a(x,y)u_x)_x + (b(x,y)u_y)_y = \sin(\pi xy),\quad (x,y)\in[0,1]^2,
	\end{align*}
	with 
	\begin{align*}
	a(x,y) = b(x,y) = \left\lbrace\begin{array}{ll}
	10^{-3}, & (x,y)\in [0,0.5]\times [0.5,1], \\
	10^3, & (x,y) \in [0.5,1]\times [0,0.5], \\
	1, & \text{otherwise},
	\end{array}\right.
	\end{align*}
	with a Dirichlet boundary condition discretized with central differences.
\end{itemize}
A regular mesh was assumed in constructing the finite difference matrices for these PDEs.

The off-the-shelf matrices are from the University of Florida Sparse Matrix Collection \cite{davis2011university}: 
we used all symmetric matrices having smaller than 65536 rows/columns. 
The matrices come from a variety of scientific problems: structural engineering, 
theoretical/quantum chemistry, heat flow, 3D vision, finite element approximations and networks. 
To enable application of the Daubechies wavelet transform for the implicit wavelet preconditioner, 
we discarded a random set of rows/columns from each matrix such that its size is reduced to 
$p\ts 2^s$, where $s = \lfloor \log_2 n\rfloor$ and $p = \lfloor n/2^s\rfloor$. 
The right hand sides of the linear systems were random vectors drawn from a multivariate normal 
distribution with mean zero and unit variance. 

For the model PDE problems, we used GMRES with a stopping tolerance of $10^{-8}$ in relative 
residual and a cap on the number of iterations at 1000. 
For the off-the-shelf matrices, we use a tolerance of $10^{-4} $ and an iterations cap of 500. 
We only show those matrices for which GMRES convergence was achieved for at least one of the employed preconditioning methods (including no preconditioning). 

We implemented both wavelet preconditioners in MATLAB and parallelized the code. Daubechies wavelets \cite{daubechies1992ten} were used for both the wavelet sparse approximate preconditioners. For the model problems, we used wavelet transforms of the same dimension as the underlying PDE (whenever applicable and whenever known) while for the off-the-shelf matrices, we used one dimensional wavelet transforms. The number of wavelet levels used was 8. 

The pMMF library \cite{corr/KondorTM15} was used to compute the MMF preconditioner. Default parameters supplied by the library were used. These include using second order rotations, i.e., Givens rotations, designating half of the active number of columns at each level as wavelets and compressing the matrix until the core is of size $100\times 100$. The parameter which controls the extent of pMMF parallelization, namely the maximum size of blocks in blocked matrices, was set to 2000.

MMF preconditioning is consistently better on model problems in terms of iteration count. 
Higher dimensional finite difference Laplacian matrices are generally well conditioned, as the condition 
number depends more strongly on the mesh size. 
In fact, the condition number of $d$-dimensional finite difference Laplacian matrix grows as $n^\frac{2}{h}$. 
Even on higher dimensional Laplacians, where the  wavelet preconditioners fail to provide adequate speedup, 
MMF preconditioning is effective. 
On average, MMF preconditioning seems to converge in about half the number of iterations as that required by the best 
wavelet preconditioner. The iteration counts are tabulated in \cref{fig:fem_results}. 

In \cref{fig:fem_timing_results}, we present the wall clock running times for linear solves with the different preconditioners. In terms of the total time for the linear solve including preconditioner setup, MMF preconditioner is consistently better. Note that we used the most basic parameters while computing the MMF. With proper tuning, performance can be brought up, which would result in better performance. The other wavelet preconditioners have only one parameter, namely the level of the wavelet transform, which leaves little room for tuning. 

Increasing the wavelet transform level increases the accuracy of the wavelet preconditioners. In this case, Hawkins and Chen \cite{hawkins2005implicit} remark that a few iterations of GMRES can be used in place of reduced QR factorization in Step 4 of \cref{alg:hawkins2} to alleviate the increased setup time. However, using GMRES defeats the purpose of maintaining higher accuracy with a higher wavelet transform level.

\input{tbl-fem_results.tex}
\input{tbl-fem_timing_results.tex}

\begin{figure}
	\label{fig:fem_plots}
\begin{centering}
	\includegraphics[width=0.49\textwidth]{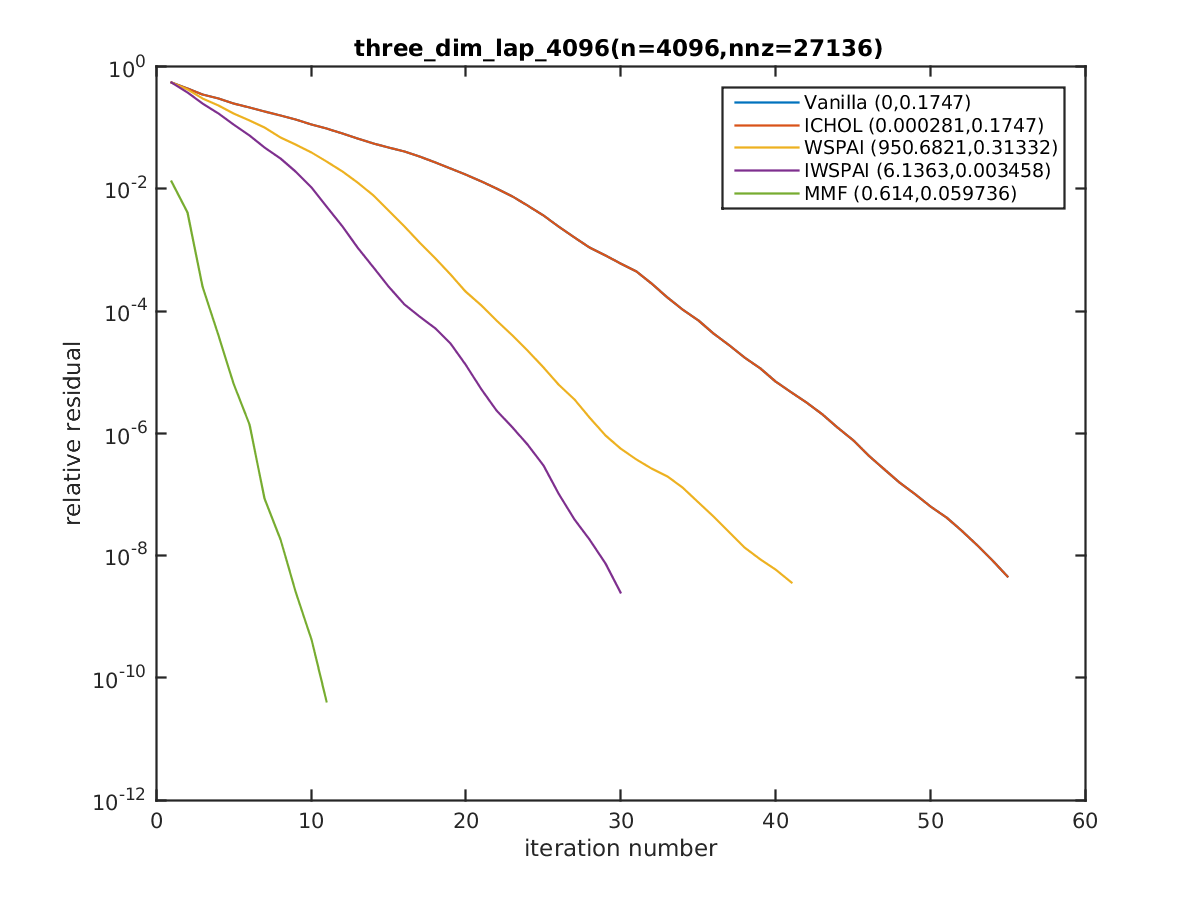}
	\includegraphics[width=0.49\textwidth]{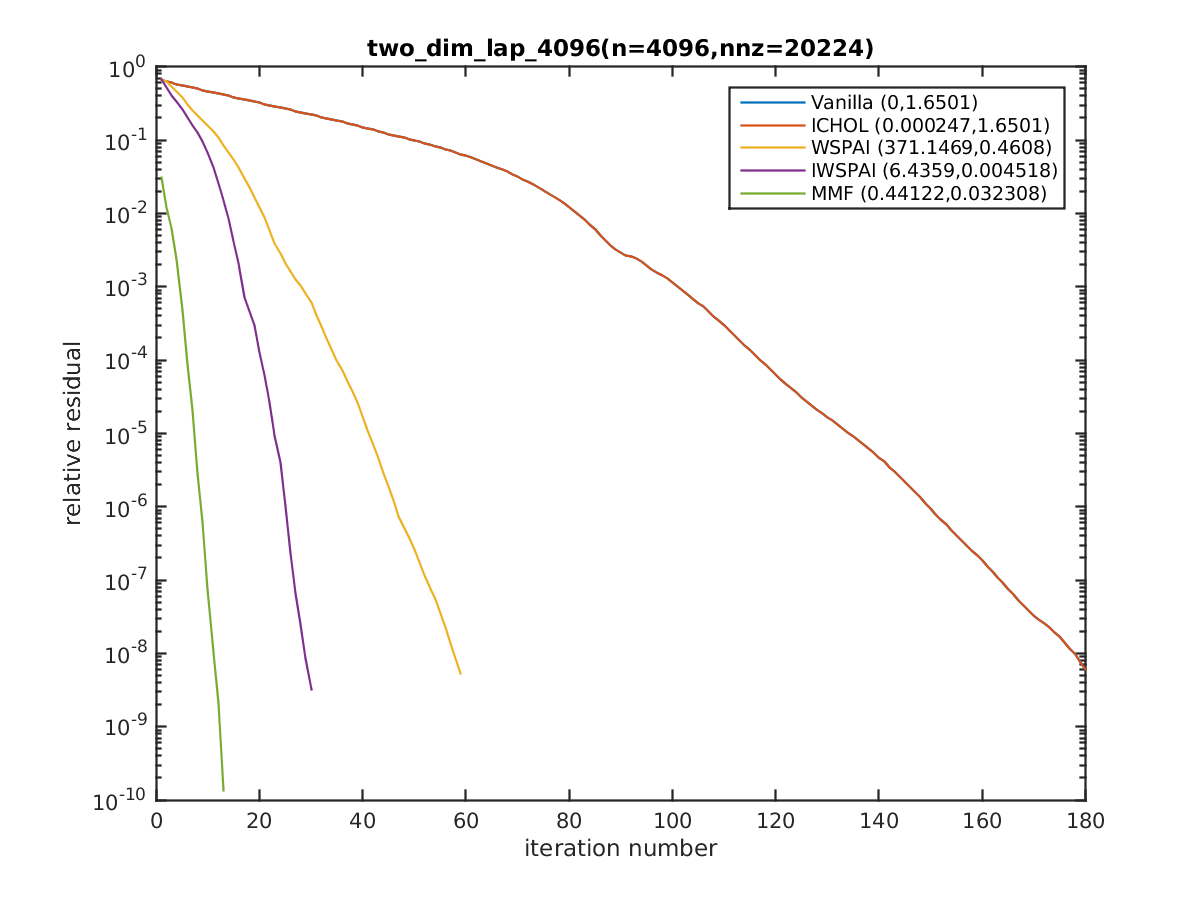}
\end{centering}
	\caption{Relative residual as a function of iteration number.}
\end{figure}

In applications where only an approximate solution to the linear system is required, it is important that 
the preconditioner lead to a reasonably accurate solution in just a small number of iterations. 
In \cref{fig:fem_plots} we plot the relative residual as a function of the iteration number. 
Relative residual is defined as $\frac{\norm{Ax_n-b}}{b}$, where $x_n$ is the $n$-th iterate. 
We see that the curve corresponding to the MMF preconditioner is below the curves for the other preconditioners. This means that an approximate solution can be determined quickly by the MMF preconditioner. 

\ignore{
In \cref{fig:qmats} we show the sparsity pattern of the $Q=Q_L Q_{L-1} \ldots Q_1$ matrices of MMF 
(see \cref{eq:mmf}) for the finite difference matrices of model PDEs 1--3. 
Hawkins and Chen \cite{hawkins2005implicit} prove that a sparse approximate inverse computed by solving \cref{eq:hawkins_chen_minimization_problem} possesses a band pattern similar to what we see in the figure (Theorem 6.1, \cite{hawkins2005implicit}). MMF naturally encodes this so-called one-sided finger pattern.

\begin{figure}
	\label{fig:qmats}
	\begin{centering}
		\includegraphics[width=0.33\textwidth]{qfigs/1D_4096.png}
		\includegraphics[width=0.32\textwidth]{qfigs/2D_4096.png}
		\includegraphics[width=0.33\textwidth]{qfigs/3D_4096.png}
	\end{centering}
	\caption{Sparsity patterns of the $Q$ matrix in the MMF of the 1D, 2D and 3D Laplacian matrices respectively of size $4096 \times 4096$.}
\end{figure}}

For the off-the-shelf matrices, we only consider the implicit wavelet preconditioner of Hawkins and Chen 
\cite{hawkins2005implicit} for comparison. 
The original wavelet preconditioner of Chan, Tang and Wan \cite{chan1997wavelet} is too slow for these large matrices. 

\input{tbl-ufl-results.tex}

In \cref{fig:ufl_results} we compare the iteration counts. The best result for each dataset is highlighted in 
bold. In the majority of datasets, the MMF preconditioner turns out best. 
However, for a few datasets such as \texttt{gyro\_m}, \texttt{crystm03}, \texttt{crystm02}, 
the implicit wavelet preconditioner outperforms MMF preconditioning. 

We remark that the ``geometry free'' nature of MMF preconditioner makes it more flexible than standard 
wavelet preconditioners. In particular, MMF can be applied to matrices of any size, not just $p\ts 2^s$. 
Furthermore, MMF preconditioning is completely invariant to the ordering of the rows/columns, in contrast 
to, for example, the multiresolution preconditioner of Bridson and Tang \cite{bridson2001multiresolution}. The adaptability of MMF makes it suitable to preconditioning a wide variety of linear systems. 

%% file: tbl-fem_results.tex
\begin{figure}
	\label{fig:fem_results}
	\centering
	\small 
	\begin{tabular}{|cc|c|c|c|c|}
		\hline
		Dataset & $n$ & no prec. & WSPAI & IWSPAI & MMF prec.\\
		\hline
		1D Laplacian & 256 		 & 256 & 46 & 13 & \textbf{10} \\
					 & 512 & 512 & 64 & 13 & \textbf{10} \\
					 & 1024	 & 1001 & 93 & 17 & \textbf{13} \\
					 & 2048  & 1001 & 131 & 17 & \textbf{2} \\
\hline
		2D Laplacian & 256  & 45 & 33 & 28 & \textbf{8} \\
					 & 1024  & 91 & 41 & 28 & \textbf{8} \\
					 & 4096  & 180 & 59 & 30 & \textbf{13} \\
\hline
		3D Laplacian & 512  & 28 & 26 & 28 & \textbf{8} \\
					 &4096  & 55 & 41 & 30 & \textbf{11} \\
\hline
		2D Disc		 & 256 &  240 & 256 & 37 & \textbf{13} \\
					 & 1024 & 868 & $\times$ & 24 & \textbf{13} \\
		\hline
	\end{tabular}
	\caption{Iteration counts of GMRES until convergence to a relative residual of $10^{-8}$. 
Here $n$ is the number of rows of the finite difference matrix. 
WSPAI refers to the wavelet sparse preconditioner of Chan, Tang and Wan \cite{chan1997wavelet} and IWSPAI to 
the implicit sparse preconditioner of Hawkins and Chen \cite{hawkins2005implicit}. 
It is clear that MMF preconditioner is consistently better. 
$\times$ indicates that the desired tolerance was not reached within 1000 iterations.}
\end{figure}

%% file: tbl-fem_timing_results.tex
\begin{figure}[t]
	\label{fig:fem_timing_results}
	\centering
	\small
	\begin{tabular}{|cc|c|cc|cc|cc|}
		\hline
		Dataset & $n$ & no prec. & \multicolumn{2}{c}{WSPAI} & \multicolumn{2}{|c|}{IWSPAI} & \multicolumn{2}{c|}{MMF prec.}  \\
		& & solve & setup & solve  & setup & solve  & setup & solve  \\
		\hline
		1D Laplacian & 256 & 0.3 & 0.77 & 0.01 & 0.8 & \textbf{2e-05} & \textbf{0.01} & 0.01 \\
		& 512 & 1.35 & 1.70 & 0.03 & 1.73 & \textbf{4.3e-05} & \textbf{0.03} & 0.02 \\
		& 1024	& 5.18 & 5.36 & 0.09  & 3.79 & \textbf{8.2e-05}  & \textbf{0.07} & 0.02 \\
		& 2048 & 7.80 & 24.2 & 0.24 & 9.9 & \textbf{1.5e-04} & \textbf{0.15} & 0.02 \\
\hline
		2D Laplacian 
		& 256 & 0.54 & 21 & 0.03 & 0.26 & \textbf{3.4e-05} & \textbf{0.05} & 0.04 \\
		& 1024 & 0.08 & 3.87 & 0.03 & 0.32 & \textbf{2.4e-04} & \textbf{0.10} & 0.02 \\
		& 4096 & 1.65 & 371 & 0.46 & 6.43 & \textbf{4.5e-03} & \textbf{0.44} & 0.03 \\
\hline 
		3D Laplacian 
		& 512 & 0.01 & 0.16 & 0.01 & 0.16 & \textbf{1.2e-05} & \textbf{0.04} & 0.01 \\
		& 4096 & 0.17 & 950 & 0.31 & 6.13 & \textbf{3.4e-03} & \textbf{0.61} & 0.05 \\
\hline 
		2D Disc 
		& 256 & 0.23 & 0.18 & 0.30 & 0.20 & \textbf{3.3e-05} & \textbf{0.01} & 0.02 \\
		& 1024 & 2.67 & 3.77 & 5.60 & 0.31 & \textbf{2.9e-04} & \textbf{0.11} & 0.03 \\
		& 4096 & 3.96 & $\times$ & $\times$ & 3.27 & 3.5e-03 & \textbf{0.41} & \textbf{2.5e-03} \\		
		\hline
	\end{tabular}
	\caption{Wall clock running time of preconditioner setup and linear solve times in seconds.  
$\times$ indicates that the desired tolerance was not reached within 1000 iterations.}
\end{figure}

%% file: tbl-ufl-results.tex
\begin{figure}
		\label{fig:ufl_results}
	\centering
	\small
	\begin{tabular}{|c|c|c|c|c|}
		\hline
		Dataset & $n$ & no prec. & IWSPAI & MMF prec.   \\
		\hline
		nd3k	&	8192	&	455	&	\textbf{236}	&	323\\
		nemeth03	&	9216	&	4	&	4	&	\textbf{2}\\
		net25	&	9216	&	\textbf{460}	&	$\times$	&	$\times$\\
		fv2	&	9216	&	\textbf{20}	&	\textbf{20}	&	27\\
		fv3	&	9216	&	42	&	\textbf{38}	&	52\\
		nemeth12	&	9216	&	13	&	10	&	\textbf{3}\\
		nemeth11	&	9216	&	10	&	8	&	\textbf{3}\\
		nemeth09	&	9216	&	7	&	6	&	\textbf{3}\\
		nemeth14	&	9216	&	$\times$	&	$\times$	&	\textbf{8}\\
		nemeth04	&	9216	&	5	&	4	&	\textbf{3}\\
		nemeth23	&	9216	&	\textbf{211}	&	$\times$	&	$\times$\\
		pf2177	&	9216	&	\textbf{174}	&	$\times$	&	$\times$\\
		bloweybq	&	9216	&	$\times$	&	\textbf{8}	&	$\times$\\
		nemeth10	&	9216	&	8	&	7	&	\textbf{3}\\
		flowmeter0	&	9216	&	$\times$	&	$\times$	&	\textbf{9}\\
		nemeth25	&	9216	&	\textbf{164}	&	$\times$	&	$\times$\\
		nemeth24	&	9216	&	\textbf{179}	&	$\times$	&	$\times$\\
		nemeth15	&	9216	&	282	&	$\times$	&	\textbf{70}\\
		nopoly	&	10240	&	119	&	108	&	\textbf{105}\\
		bcsstk17	&	10240	&	$\times$	&	$\times$	&	\textbf{266}\\
		bundle1	&	10240	&	$\times$	&	$\times$	&	\textbf{30}\\
		linverse	&	11264	&	$\times$	&	\textbf{20}	&	$\times$\\
		t2dah	&	11264	&	$\times$	&	$\times$	&	\textbf{7}\\
		crystm02	&	13312	&	\textbf{1}	&	\textbf{1}	&	30\\
		Pres\_Poisson	&	14336	&	436	&	\textbf{43}	&	114\\
		bcsstm25	&	14336	&	$\times$	&	$\times$	&	\textbf{2}\\
		gyro\_m	&	16384	&	\textbf{1}	&	\textbf{1}	&	115\\
		gyro\_k	&	16384	&	$\times$	&	$\times$	&	\textbf{220}\\
		nd6k	&	16384	&	$\times$	&	\textbf{270}	&	330\\
		bodyy4	&	16384	&	184	&	147	&	\textbf{91}\\
		t3dl\_a	&	18432	&	$\times$	&	141	&	\textbf{6}\\
		Si5H12	&	18432	&	103	&	\textbf{71}	&	89\\
		Trefethen\_20000b	&	18432	&	$\times$	&	$\times$	&	\textbf{8}\\
		crystm03	&	24576	&	\textbf{1}	&	\textbf{1}	&	33\\
		spmsrtls	&	28672	&	$\times$	&	\textbf{150}	&	$\times$\\
		wathen100	&	28672	&	$\times$	&	$\times$	&	\textbf{33}\\
		wathen120	&	32768	&	$\times$	&	$\times$	&	\textbf{33}\\
		mario001	&	36864	&	269	&	$\times$	&	$\times$\\
		torsion1	&	36864	&	41	&	\textbf{29}	&	50\\
		bfly	&	49152	&	\textbf{59}	&	$\times$	&	$\times$\\
		crankseg\_2	&	57344	&	$\times$	&	$\times$	&	\textbf{246}\\
		Ga3As3H12	&	57344	&	$\times$	&	$\times$	&	\textbf{104}\\
		cant	&	57344	&	$\times$	&	$\times$	&	\textbf{83}\\   
		\hline
	\end{tabular}
	\caption{Iteration counts of GMRES solved to a relative error of $10^{-4}$. 
$\times$ indicates that the method did not achieve the desired tolerance within 500 iterations.}
\end{figure}

%% file: conclusion.tex
\section{Conclusion}
\label{sec:conclusion}

We presented a new multiresolution preconditioner for symmetric linear systems that does not depend 
on any geometric assumptions, and hence can be applied to any coefficient matrix 
that is assumed to have multiresolution structure, even in the loose sense. 
Numerical experiments show the effectiveness of the new preconditioner in a range of problems. 
In our experiments we used default parameters, but with fine tuning our results could possibly be 
improved further. 

It is not yet clear exactly what kind of matrices the new MMF preconditioner is most 
effective on, in part due to the general nature of the pMMF algorithm.  It is possible that specializing MMF to specific types of linear systems would yield even more 
effective preconditioners. 

\section{Acknowledgements}
We would like to thank Prof.\;Stuart Hawkins for help with implementing his preconditioner and Prof.\;Jonathan Weare 
for discussions. This work was funded by NSF award CCF--1320344.